\documentclass{article}
\usepackage[utf8]{inputenc}
\usepackage{fullpage}
\usepackage{amsthm,mathtools,scalerel}
\usepackage{color,soul,latexsym,amsmath,amssymb,amsfonts,amsthm,dsfont,enumitem,xcolor,bbm,threeparttable,graphicx,float,subfigure}
\usepackage{authblk}
\usepackage[round]{natbib}
\RequirePackage[colorlinks=true,allcolors=blue,hypertexnames=false]{hyperref}%
\usepackage{graphicx,todonotes}
\usepackage{outlines}
\usepackage{array}
\usepackage{selectp}

\usepackage[ruled,linesnumbered, vlined, noend]{algorithm2e}
\mathtoolsset{showonlyrefs}
\usepackage{titlesec}

\newtheorem{thm}{Theorem}
\newtheorem{prop}{Proposition}

\DeclareMathOperator*{\argmin}{argmin}

\title{\bf Median bias of M-estimators}
\author{
    Arun Kumar Kuchibhotla}
\affil{{\tt arunku@cmu.edu}}
\affil{Department of Statistics \& Data Science, Carnegie Mellon University}
\date{}                     

\begin{document}

\maketitle

\begin{abstract}
In this note, we derive bounds on the median bias of univariate $M$-estimators under mild regularity conditions. These requirements are not sufficient to imply convergence in distribution of the M-estimators. We also discuss median bias of some multivariate M-estimators.
\end{abstract}
\maketitle

\section{Introduction}
Commonly used estimators, known as M-estimators, are obtained as solution to optimization problems. Under certain regularity conditions, properly normalized M-estimators are shown to convergence in distribution to a mean zero Gaussian. This implies that the median of the M-estimator converges to the population parameter, i.e., the M-estimator is asymptotically median unbiased. The aim of this note is to show that in some cases asymptotic median unbiasedness can be proved without proving convergence in distribution. The more interesting aspect is that median unbiasedness can be proved under far less regularity conditions than those required for convergence in distribution.

For any estimator $\widehat{\theta}$ estimating $\theta_0$, we define the median bias as
\[
\mbox{Med-bias}_{\theta_0}(\widehat{\theta}) := \left(\frac{1}{2} - \max\left\{\mathbb{P}(\widehat{\theta} \le \theta_0),\,\mathbb{P}(\widehat{\theta} \ge \theta_0)\right\}\right)_+,
\]
where $(x)_+ := \max\{x, 0\}$.
If $\mathbb{P}(\widehat{\theta} \ge \theta_0) \ge 1/2$ and $\mathbb{P}(\widehat{\theta} \le \theta_0) \ge 1/2$, then $\mbox{Med-bias}_{\theta_0}(\widehat{\theta}) = 0$ and $\widehat{\theta}$ is median unbiased for $\theta_0$.
\paragraph{Notation.} For any function $f$, we use $\dot{f}$ and $\ddot{f}$ to denote the first and second derivatives of $f$. 
\section{Univariate M-estimators}
In this section, we consider median bias of univariate M-estimators and Z-estimators. 
\subsection{Convex M-estimators}
Suppose $\Theta\subseteq\mathbb{R}$ and $\mathbb{M}_n:\Theta\to\mathbb{R}$ is a convex function and define
\[
\widehat{\theta}_n ~:=~ \argmin_{\theta\in\Theta}\,\mathbb{M}_n(\theta).
\]
If $\mathbb{M}_n(\cdot)$ converges in probability pointwise to $M(\cdot)$, then we can define the target $\theta_0$ of $\widehat{\theta}_n$ as the minimizer of $M(\theta)$. Formally, $\theta_0 = \argmin_{\theta\in\Theta}M(\theta)$. With $\dot{\mathbb{M}}_n(\cdot)$ representing the derivative of $\mathbb{M}_n(\cdot)$, convexity of $\mathbb{M}_n(\cdot)$ and an assumption that $\theta_0$ lies in the interior of $\Theta$ implies that
\begin{equation}\label{eq:overestimation-inequality}
\dot{\mathbb{M}}_n(\theta_0) < 0 \quad\Rightarrow\quad \widehat{\theta}_n \ge \theta_0 \quad\Rightarrow\quad \dot{\mathbb{M}}_n(\theta_0) \le 0,
\end{equation}
and
\begin{equation}\label{eq:underestimation-inequality}
\dot{\mathbb{M}}_n(\theta_0) > 0 \quad\Rightarrow\quad \widehat{\theta}_n \le \theta_0 \quad\Rightarrow\quad \dot{\mathbb{M}}_n(\theta_0) \ge 0.
\end{equation}
See, for example, (3.4) of~\cite{bentkus1997berry} for the proof.
Hence, we get that 
\[
\mathbb{P}(\widehat{\theta}_n \ge \theta_0) ~\ge~ \mathbb{P}(\dot{\mathbb{M}}_n(\theta_0) < 0)\quad\mbox{and}\quad \mathbb{P}(\widehat{\theta}_n \le \theta_0) ~\ge~ \mathbb{P}(\dot{\mathbb{M}}_n(\theta_0) > 0).
\]
These inequalities imply the following result.
\begin{thm}\label{thm:median-bias-convex}
If $\theta_0$ lies in the interior of $\Theta$ and $\theta\mapsto\mathbb{M}_n(\theta)$ is convex, then
\begin{equation}\label{eq:median-bias-bound-convex}
\mathrm{Med\mbox{-}bias}_{\theta_0}(\widehat{\theta}_n) ~\le~ \left(\frac{1}{2} - \max\left\{\mathbb{P}(\dot{\mathbb{M}}_n(\theta_0) < 0),\,\mathbb{P}(\dot{\mathbb{M}}_n(\theta_0) > 0)\right\}\right)_+.
\end{equation}
\end{thm}
Theorem~\ref{thm:median-bias-convex} implies that the median bias of $\widehat{\theta}_n$ can be controlled by studying $\dot{\mathbb{M}}_n(\theta_0)$. The study of $\dot{\mathbb{M}}_n(\theta_0)$ is not enough to prove consistency or convergence in distribution of $\widehat{\theta}_n$. Proving consistency requires assumptions on the curvature of $\mathbb{M}_n(\cdot)$ around $\theta_0$ and similarly proving asymptotic normality requires assumptions on the first derivative of $\mathbb{M}_n(\cdot)$.

In many cases, $\mathbb{M}_n(\cdot)$ is an average of random variables and hence, $\dot{\mathbb{M}}_n(\theta_0)$ satisfies a central limit theorem under Lindeberg type conditions. This implies that the right hand side of~\eqref{eq:median-bias-bound-convex} converges to zero and proves that $\widehat{\theta}_n$ is asymptotically median unbiased. Note that it is relatively straightforward to obtain a finite sample bound on the median bias using~\eqref{eq:median-bias-bound-convex} and the Berry--Esseen bounds for sum of independent or weakly dependent random variables. For examples of such Berry--Esseen bounds, see~\citet[Chap. V]{petrov2012sums} and~\cite{hormann2009berry}. Below, we provide some simple applications of Theorem~\ref{thm:median-bias-convex}.
\paragraph{Median Estimation.} The sample median of $X_1, \ldots, X_n$ is defined as the minimizer of
\[
\mathbb{M}_n(\theta) ~:=~ \sum_{i=1}^n |X_i - \theta|\quad\Rightarrow\quad \dot{\mathbb{M}}_n(\theta) = \sum_{i=1}^n\left\{2\mathbbm{1}\{X_i \le \theta\} - 1\right\}.
\]
If $X_1, \ldots, X_n$ are independent but possibly non-identically distributed observations and $\theta_0$ is a solution of $\mathbb{E}[\dot{\mathbb{M}}_n(\theta)] = 0$, then $n^{-1/2}\dot{\mathbb{M}}_n(\theta_0)$ converges in distribution to a mean zero normal distribution. Proving convergence in distribution of the sample median requires an assumption on the H{\"o}lder continuity of the distribution functions. See, for example,~\cite{knight1998limiting} and~\cite{knight1999asymptotics}. Note that the calculations above also apply to quantile estimation and proves asymptotic median unbiasedness.
\paragraph{$L_p$-median.} Generalizing the sample median, consider $\widehat{\theta}_n$ as a minimizer of $\mathbb{M}_n(\theta)$ over $\theta\in\mathbb{R}$ where for $p\ge 1$,
\[
\mathbb{M}_n(\theta) ~:=~ \sum_{i=1}^n |X_i - \theta|^{p} \quad\Rightarrow\quad \dot{\mathbb{M}}_n(\theta) = p\sum_{i=1}^n |X_i - \theta|^{p-1}\mbox{sign}(\theta - X_i).
\]
Define $\theta_0$ as a solution to the equation $\mathbb{E}[\dot{\mathbb{M}}_n(\theta)] = 0$. Inequality~\eqref{eq:median-bias-bound-convex} along with the Berry--Esseen bounds shows that $\widehat{\theta}_n$ is asymptotically median unbiased for $\theta_0$. Once again more conditions on the distribution of $X_i$'s is needed (for $p\le 3$) to ensure convergence in distribution; see~\cite{bentkus1997berry}.
\paragraph{Maximum Likelihood Estimator.} Suppose $\{p_{\theta}:\,\theta\in\Theta\}$ is a family of parametric densities parametrized by $\theta\in\Theta\subseteq\mathbb{R}$. Consider the maximum likelihood estimator (MLE) $\widehat{\theta}_n$ as
\[
\widehat{\theta}_n ~:=~ \argmin_{\theta\in\Theta}\,-\sum_{i=1}^n \log p_{\theta}(X_i).
\]
If $\theta\mapsto-\log p_{\theta}(x)$ is convex, then $\mathbb{M}_n(\theta) = -\sum_{i=1}^n \log p_{\theta}(X_i)$ is a convex function of $\theta$. Assuming differentiability in quadratic mean (DQM) of the parametric family, let $u_{\theta}(x)$ be the likelihood score function; see Eq. (7.1) of~\cite{van2000asymptotic} for DQM. Define $\theta_0$ to be a solution to the equation $\sum_{i=1}^n \mathbb{E}[u_{\theta}(X_i)] = 0$. Assuming $X_1, \ldots, X_n$ are independent and the Linderberg condition on $u_{\theta_0}(X_i), 1\le i\le n$, we get that the MLE is asymptotically median unbiased. Once again a properly normalized MLE need not converge in distribution without further assumptions that guarantee asymptotic equicontinuity of the likelihood score. Further under possible misspecification, the Jacobian also needs to be non-zero at $\theta_0$ for convergence in distribution.
\subsection{Non-differentiable M-estimators}
We have assumed the existence of a version of the derivative $\dot{\mathbb{M}}_n(\cdot)$ in the previous subsections. It is possible to avoid such assumption. From the convexity of $\mathbb{M}_n(\cdot)$, it follows that for any $\varepsilon > 0$,
\[
\{\widehat{\theta}_n > \theta_0 + \varepsilon\}\quad\Rightarrow\quad \{\mathbb{M}_n(\theta_0) \ge \mathbb{M}_n(\theta_0 + \varepsilon)\},
\]
and
\[
\{\widehat{\theta}_n < \theta_0 - \varepsilon\}\quad\Rightarrow\quad \{\mathbb{M}_n(\theta_0) \ge \mathbb{M}_n(\theta_0 - \varepsilon)\}.
\]
Hence,
\begin{align*}
\mathbb{P}(\widehat{\theta}_n \le \theta_0 + \varepsilon) ~&\ge~ \mathbb{P}(\mathbb{M}_n(\theta_0) < \mathbb{M}_n(\theta_0 + \varepsilon)),\\
\mathbb{P}(\widehat{\theta}_n \ge \theta_0 - \varepsilon) ~&\ge~ \mathbb{P}(\mathbb{M}_n(\theta_0) < \mathbb{M}_n(\theta_0 - \varepsilon)).
\end{align*}
Now observe that $$\mathbb{P}(\widehat{\theta}_n \le \theta_0) = \lim_{\varepsilon\downarrow0}\mathbb{P}(\widehat{\theta}_n \le \theta_0 + \varepsilon)\quad\mbox{and}\quad \mathbb{P}(\widehat{\theta}_n \ge \theta_0) = \lim_{\varepsilon\downarrow0}\mathbb{P}(\widehat{\theta}_n \ge \theta_0 - \varepsilon).$$
Therefore,
\begin{equation}\label{eq:median-bias-non-differentiable}
\mbox{Med-bias}_{\theta_0}(\widehat{\theta}) ~\le~ \lim_{\varepsilon\downarrow0}\left(\frac{1}{2} - \max\left\{\mathbb{P}(\mathbb{M}_n(\theta_0) < \mathbb{M}_n(\theta_0 + \varepsilon)),\,\mathbb{P}(\mathbb{M}_n(\theta_0) < \mathbb{M}_n(\theta_0 - \varepsilon))\right\}\right)_+.
\end{equation}
In many cases, $\mathbb{M}_n(\theta)$ is an average of random variables and $\mathbb{M}_n(\theta_0) - \mathbb{M}_n(\theta_0 + \varepsilon)$ converges to a negative quantity for any fixed $\varepsilon > 0$.
\paragraph{Maximum Likelihood Estimation.}
Consider again the problem of median bias of the maximum likelihood estimator. In this case, $\mathbb{M}_n(\theta) = -\sum_{i=1}^n \log p_{\theta}(X_i)$ and the target of the MLE $\theta_0$ is defined as the minimizer of $\theta\mapsto\mathbb{M}_n(\theta)$ over $\theta\in\Theta$. If $\theta_0$ lies in the interior of $\Theta$ and $\theta_0\pm\varepsilon\in\Theta$, then
\begin{align*}
&\mathbb{P}(\mathbb{M}_n(\theta_0) < \mathbb{M}_n(\theta_0 + \varepsilon))\\ 
&\quad= \mathbb{P}\left(\sum_{i=1}^n \log\frac{p_{\theta_0 + \varepsilon}(X_i)}{p_{\theta_0}(X_i)} < 0\right)\\
&\quad= \mathbb{P}\left(\sum_{i=1}^n \left\{\log\frac{p_{\theta_0 + \varepsilon}(X_i)}{p_{\theta_0}(X_i)} - \mathbb{E}\left[\log\frac{p_{\theta_0 + \varepsilon}(X_i)}{p_{\theta_0}(X_i)}\right]\right\} < -\sum_{i=1}^n\mathbb{E}\left[\log\frac{p_{\theta_0 + \varepsilon}(X_i)}{p_{\theta_0}(X_i)}\right]\right).
\end{align*}
Note that, by definition, $\sum_{i=1}^n \mathbb{E}\log(p_{\theta_0 + \varepsilon}(X_i)/p_{\theta_0}(X_i)) \le 0$ and is strictly negative if $\varepsilon > 0$ by identifiability. Hence, it follows that
\[
\mathbb{P}(\mathbb{M}_n(\theta_0) < \mathbb{M}_n(\theta_0 + \varepsilon)) \ge \mathbb{P}\left(\sum_{i=1}^n \left\{\log\frac{p_{\theta_0 + \varepsilon}(X_i)}{p_{\theta_0}(X_i)} - \mathbb{E}\left[\log\frac{p_{\theta_0 + \varepsilon}(X_i)}{p_{\theta_0}(X_i)}\right]\right\} \le 0\right).
\]
Similarly,
\[
\mathbb{P}(\mathbb{M}_n(\theta_0) < \mathbb{M}_n(\theta_0 - \varepsilon)) \ge \mathbb{P}\left(\sum_{i=1}^n \left\{\log\frac{p_{\theta_0 - \varepsilon}(X_i)}{p_{\theta_0}(X_i)} - \mathbb{E}\left[\log\frac{p_{\theta_0 - \varepsilon}(X_i)}{p_{\theta_0}(X_i)}\right]\right\} \le 0\right).
\]
If the log-likelihood ratio satisfies the Linderbeg condition, then~\eqref{eq:median-bias-non-differentiable} implies that the MLE is again asymptotically median unbiased.
\subsection{Non-convex M-estimators}
Convexity of the objective function $\mathbb{M}_n(\cdot)$ is not very crucial for~\eqref{eq:median-bias-bound-convex}. All that is required is that $\theta\mapsto\mathbb{M}_n(\theta)$ is convex in the neighborhood of $\theta_0$ and with some positive probability the estimator $\widehat{\theta}_n$ belongs to that neighborhood. These conditions are same as convexity and consistency assumptions (1.4), (1.5) in~\cite{bentkus1997berry}. Formally, for set
\begin{align*}
    \eta_{1,n}(\delta) ~&:=~ 1 - \mathbb{P}\left(\theta\mapsto\mathbb{M}_n(\theta)\mbox{ is convex on }[\theta_0 - \delta, \theta_0 + \delta]\right),\\
    \eta_{2n}(\delta) ~&:=~ \mathbb{P}(|\widehat{\theta}_n - \theta_0| > \delta).
\end{align*}
On the event $\widehat{\theta}_n\in[\theta_0 - \delta, \theta_0 + \delta]\subseteq\Theta$, inequalities~\eqref{eq:overestimation-inequality} and~\eqref{eq:underestimation-inequality} hold true. Therefore, we get
\begin{equation}\label{eq:median-bias-bound-non-convex}
    \mbox{Med-bias}_{\theta_0}(\widehat{\theta}_n) ~\le~ \left(\frac{1}{2} - \max\left\{\mathbb{P}(\dot{\mathbb{M}}_n(\theta_0) < 0),\,\mathbb{P}(\dot{\mathbb{M}}_n(\theta_0) > 0)\right\}\right)_+ + \min_{\delta\ge0}\left[\eta_{1,n}(\delta) + \eta_{2,n}(\delta)\right]. 
\end{equation}
Note that if $\theta\mapsto\mathbb{M}_n(\theta)$ is convex on $\Theta$, then $\eta_{1,n}(\infty) = 0$ and $\eta_{2,n}(\infty) = 0$. Inequality~\eqref{eq:median-bias-bound-non-convex} follows from (3.3) of~\cite{bentkus1997berry}.

Alternatively, one can consider the usual Taylor series expansion way and prove a better result. For this, we additionally require absolute continuity of the first derivative of $\mathbb{M}_n(\cdot)$. If $\widehat{\theta}_n$ solves the equation $\dot{\mathbb{M}}_n(\theta) = 0$ and using absolute continuity of $\theta\mapsto\dot{\mathbb{M}}_n(\theta)$, we get that
\[
0 = \dot{\mathbb{M}}_n(\theta_0) + \int_0^1 \ddot{\mathbb{M}}_n(\theta_0 + t(\widehat{\theta}_n - \theta_0))dt (\widehat{\theta}_n - \theta_0).
\]
Now assuming that $\widehat{\theta}_n$ is a locally unique solution with probability 1, we conclude that $\int_0^1 \ddot{\mathbb{M}}_n(\theta_0 + t(\widehat{\theta}_n - \theta_0))dt \neq 0$ and hence,
\begin{equation}\label{eq:Taylor-series-M-estimator}
\widehat{\theta}_n - \theta_0 = -\frac{\dot{\mathbb{M}}_n(\theta_0)}{\int_0^1 \ddot{\mathbb{M}}_n(\theta_0 + t(\widehat{\theta}_n - \theta_0))dt}.    
\end{equation}
This implies that 
\begin{equation}\label{eq:median-bias-Z-estimator}
\mbox{Med-bias}_{\theta_0}(\widehat{\theta}_n) ~=~ \left(\frac{1}{2} - \max\left\{\mathbb{P}(\dot{\mathbb{M}}_n(\theta_0) \ge 0),\,\mathbb{P}(\dot{\mathbb{M}}_n(\theta_0) \le 0)\right\}\right)_+.
\end{equation}
Equation~\eqref{eq:Taylor-series-M-estimator} provides a more intuitive reason for why median bias of \textit{univariate} M-estimators can be controlled without requiring conditions to imply convergence in distribution. To prove convergence in distribution, we would require $\widehat{\theta}_n$ to be consistent for $\theta_0$ along with conditions to ensure that the denominator on the right hand side of~\eqref{eq:Taylor-series-M-estimator} can be replaced with $\ddot{\mathbb{M}}_n(\theta_0)$. 

The bound on median bias~\eqref{eq:median-bias-Z-estimator} readily applies to Z-estimators which are obtained as solutions of estimating equations rather than minimizers of objective functions.
\section{Multivariate M-estimation}
The calculations from previous sections can be used trivially when estimating a parameter in presence of nuisance parameters. Suppose $\mathbb{M}_n:\Theta\times\Lambda\to\mathbb{R}$ be an objective function and define
\[
(\widehat{\theta}_n, \widehat{\lambda}_n) ~:=~ \argmin_{\theta,\lambda}\,\mathbb{M}_n(\theta, \lambda).
\]
Setting $M(\theta,\lambda)$ as the pointwise limit in probability of $\mathbb{M}_n(\theta, \lambda)$, define
\[
(\theta_0, \lambda_0) ~:=~ \argmin_{\theta, \lambda}\,M(\theta, \lambda).
\]
To apply the results from previous section to study the median bias of $\widehat{\theta}_n$, observe that
\[
\widehat{\theta}_n ~:=~ \argmin_{\theta\in\Theta}\,\min_{\lambda}\mathbb{M}_n(\theta, \lambda).
\]
If $(\theta, \eta)\mapsto\mathbb{M}_n(\theta, \eta)$ is a convex function, then $\theta\mapsto\min_{\lambda}\,\mathbb{M}_n(\theta, \lambda)$ is also a convex function. This is called the inf-projection of $\mathbb{M}_n(\theta, \eta)$. However, $\mathbb{W}_n(\theta) = \min_{\lambda}\mathbb{M}_n(\theta, \lambda)$ and its derivative (sub-gradient) $\dot{\mathbb{W}}_n(\theta)$ are complicated functions to study, in general. In some special cases, $\mathbb{W}_n(\theta)$ is available in closed form and might be easily analysed.
\subsection{Application 1: Least Squares Linear Regression} 
Suppose $(Y_i, T_i, X_i)\in\mathbb{R}\times\mathbb{R}\times\mathbb{R}^d$, $1\le i\le n$ represent the set of observations on the treatment variable ($T$), covariates ($X$), and the response ($Y$). Consider
\begin{equation}\label{eq:OLS-regression}
(\widehat{\theta}_n, \widehat{\lambda}_n) ~:=~ \argmin_{\theta,\lambda}\,\sum_{i=1}^n (Y_i - \theta T_i - \lambda^{\top}X_i)^2.    
\end{equation}
The targets of $\widehat{\theta}_n, \widehat{\lambda}_n$ are defined as
\[
(\theta_0, \lambda_0) ~:=~ \argmin_{\theta,\lambda}\sum_{i=1}^n \mathbb{E}[(Y_i - \theta T_i - \lambda^{\top}X_i)^2].
\]
This can be written as
\[
\widehat{\theta}_n ~:=~ \argmin_{\theta}\,\sum_{i=1}^n (Y_i - \widehat{\beta}_{Y,n}^{\top}X_i - \theta(T_i - \widehat{\beta}_{T,n}^{\top}X_i))^2,
\]
where
\[
\widehat{\beta}_{Y,n} ~:=~ \argmin_{\beta\in\mathbb{R}^d}\sum_{i=1}^n (Y_i - \beta^{\top}X_i)^2,\quad\mbox{and}\quad \widehat{\beta}_{T,n} ~:=~ \argmin_{\beta\in\mathbb{R}^d}\,\sum_{i=1}^n (T_i - \beta^{\top}X_i)^2.
\]
Hence $\widehat{\theta}_n$ is the minimizer of a quadratic convex objective function and~\eqref{eq:median-bias-bound-convex} (or~\eqref{eq:median-bias-bound-non-convex}) yields a bound on the median bias of $\widehat{\theta}_n$. For notational convenience, set $\widehat{R}_{Y,i} = Y_i - \widehat{\beta}_{Y,n}^{\top}X_i$ and $\widehat{R}_{T,i} = T_i - \widehat{\beta}_{T,n}^{\top}X_i$. Also, let $\beta_{Y}$ and $\beta_T$ denote the targets of $\widehat{\beta}_{Y,n}$ and $\widehat{\beta}_{T,n}$, and set $R_{Y,i} = Y_i - \beta_{Y}^{\top}X_i$ and $R_{T,i} = R_i - \beta_{T}^{\top}X_i$.  Then~\eqref{eq:median-bias-bound-non-convex} yields
\begin{equation}\label{eq:median-bias-OLS}
\mbox{Med-bias}_{\theta_0}(\widehat{\theta}_n) = \left(\frac{1}{2} - \max\left\{\mathbb{P}\left(\sum_{i=1}^n \widehat{R}_{T,i}(\widehat{R}_{Y,i} - \theta_0\widehat{R}_{T,i}) \le 0\right),\,\mathbb{P}\left(\sum_{i=1}^n \widehat{R}_{T,i}(\widehat{R}_{Y,i} - \theta_0\widehat{R}_{T,i}) \ge 0\right)\right\}\right)_+.    
\end{equation}
This equality holds true if $T_i$ is not perfectly collinear with $X_i$. Note that unlike the examples discussed in previous sections, the sums on the right hand side of~\eqref{eq:median-bias-OLS} are not of independent random variables. 
Observe that
\begin{align*}
\sum_{i=1}^n \widehat{R}_{T,i}(\widehat{R}_{Y,i} - \theta_0\widehat{R}_{T,i}) &= \sum_{i=1}^n R_{T,i}(R_{Y,i} - \theta_0R_{T,i}) + \sum_{i=1}^n (\widehat{R}_{T,i} - R_{T,i})(R_{Y,i} - \theta_0R_{T,i})\\
&\quad- \sum_{i=1}^n \widehat{R}_{T,i}(\widehat{R}_{Y,i} - R_{Y,i} - \theta_0(\widehat{R}_{T,i} - R_{T,i})).
\end{align*}
From the definition of $\widehat{\beta}_{T,n}$ it follows that $\sum_{i=1}^n \widehat{R}_{T,i}X_i = 0$. Noting that $\widehat{R}_{Y,i} - R_{Y,i} = X_i^{\top}(\widehat{\beta}_{Y,n} - \beta_Y)$ and $\widehat{R}_{T,i} - R_{T,i} = X_i^{\top}(\widehat{\beta}_{T,n} - \beta_T)$, we conclude that
\[
\sum_{i=1}^n \widehat{R}_{T,i}(\widehat{R}_{Y,i} - R_{Y,i} - \theta_0(\widehat{R}_{T,i} - R_{T,i})) = 0.
\]
We obtain that
\[
\sum_{i=1}^n \widehat{R}_{T,i}(\widehat{R}_{Y,i} - \theta_0\widehat{R}_{T,i}) = \sum_{i=1}^n R_{T,i}(R_{Y,i} - \theta_0R_{T,i}) + (\widehat{\beta}_{T,n} - \beta_T)^{\top}\sum_{i=1}^n X_i(R_{Y,i} - \theta_0 R_{T,i}).
\]
From the definition of $\theta_0$, it can be easily verified that $\sum_{i=1}^n \mathbb{E}[R_{T,i}(R_{Y,i} - \theta_0 R_{T,i})] = 0$ and from the definitions of $\beta_{Y,n}, \beta_{T,n}$ that $\sum_{i=1}^n \mathbb{E}[X_i(R_{Y,i} - \theta_0 R_{T,i})] = 0$. Consider the event
\[
\mathcal{E}_{\eta} ~:=~ \left\{\left|(\widehat{\beta}_{T,n} - \beta_T)^{\top}\sum_{i=1}^n X_i(R_{Y,i} - \theta_0R_{T,i})\right| \le \eta\right\},
\]
and the sum
\[
S_n := \sum_{i=1}^n R_{T,i}(R_{Y,i} - \theta_0 R_{T,i}).
\]
Using this event and inequality~\eqref{eq:median-bias-OLS}, we obtain the following result.
\begin{prop}
For any sequence of random vectors $(Y_i, T_i, X_i)\in\mathbb{R}\times\mathbb{R}\times\mathbb{R}^d$, the estimator $\widehat{\theta}_n$ defined by~\eqref{eq:OLS-regression} satisfies
\[
\mathrm{Med\mbox{-}bias}_{\theta_0}(\widehat{\theta}_n) ~\le~ \inf_{\eta > 0}\left[\left(\frac{1}{2} - \max\left\{\mathbb{P}\left(S_n \le -\eta\right),\,\mathbb{P}\left(S_n \ge \eta\right)\right\}\right)_+ + \mathbb{P}(\mathcal{E}_{\eta}^c)\right].
\]
\end{prop}
Under mild moment conditions as well as weak dependence assumptions, it can be proved that $\mathbb{P}(\mathcal{E}_{\eta})$ converges to 1 as $n\to\infty$ for $\eta = O(d)$. With $\eta = O(d)$, it suffices for $d = o(\sqrt{n})$ to ensure that the median bias converges to zero. If $d = O(\sqrt{n})$, then the median bias converges to a constant bounded away from $0$ and $1/2$. The calculations above can also be used with the Neyman orthogonal estimating equation in a partial linear model~\citep{chernozhukov2018double}.
\subsection{Application 2: Partial Linear Regression with Sample Splitting} 
Suppose $(Y_i, T_i, X_i)\in\mathbb{R}\times\mathbb{R}\times\mathbb{R}^d, 1\le i\le n$ satisfy the partial linear model
\[
Y_i = \theta_0T_i + g_0(X_i) + U_i, \quad\mbox{and}\quad T_i = m_0(X_i) + V_i,\quad\mbox{where}\quad \mathbb{E}[U_i|T_i, X_i] = 0, \mathbb{E}[V_i|X_i] = 0.
\]
The parameter of interest is still $\theta_0$ (the treatment effect). Split the data into two parts $\mathcal{D}_1$ and $\mathcal{D}_2$. Let $\widehat{m}(\cdot)$ and $\widehat{g}(\cdot)$ be estimators of $m_0(\cdot)$ and $g_0(\cdot)$, respectively. We will assume that these estimators are computed from $\mathcal{D}_1$. For $i\in\mathcal{D}_2$, set $\widehat{R}_{T,i} = T_i - \widehat{m}(X_i)$ and $\widehat{R}_{Y,i} = Y_i - \widehat{g}(X_i)$. Then consider the estimator $\widehat{\theta}_n$ as a solution to the equation
\[
\mathbb{Z}_n(\theta) := \sum_{i\in\mathcal{D}_2} \widehat{R}_{T,i}(\widehat{R}_{Y,i} - D_i\theta) = 0,
\]
where $\mathcal{D}_2$ is the second part of the data. From~\eqref{eq:median-bias-Z-estimator}, it follows that
\begin{equation}\label{eq:median-bias-partial-linear}
\mbox{Med-bias}_{\theta_0}(\widehat{\theta}_n) = \left(\frac{1}{2} - \max\left\{\mathbb{P}\left(\mathbb{Z}_n(\theta_0) \le 0\right),\, \mathbb{P}(\mathbb{Z}_n(\theta_0) \ge 0)\right\}\right)_+. 
\end{equation}
Now set $\mathbb{Z}_n^c(\theta_0) = \mathbb{Z}_n(\theta_0) - \mathbb{E}[\mathbb{Z}_n(\theta_0)|\mathcal{D}_1]$.
Then,~\eqref{eq:median-bias-partial-linear} implies
\begin{equation}\label{eq:median-bias-partial-regression}
\mbox{Med-bias}_{\theta_0}(\widehat{\theta}_n) ~\le~ \left(\frac{1}{2} - \max\bigg\{\mathbb{P}(\mathbb{Z}_n^c(\theta_0) \le -|\mathbb{E}[\mathbb{Z}_n(\theta_0)|\mathcal{D}_1]|),\,\mathbb{P}(\mathbb{Z}_n^c(\theta_0) \ge |\mathbb{E}[\mathbb{Z}_n(\theta_0)|\mathcal{D}_1]|)\bigg\}\right)_+.
\end{equation}
Note that conditional $\mathcal{D}_1$ (the first split of the data), $\mathbb{Z}_n^c(\theta_0)$ is a sum of centered (mean zero) random variables. We can write
\begin{align*}
\mathbb{Z}_n(\theta_0) &= \sum_{i\in\mathcal{D}_2} R_{T,i}(R_{Y,i} - D_i\theta_0) + \sum_{i\in\mathcal{D}_2} R_{T,i}(\widehat{R}_{Y,i} - R_{Y,i})\\
&\quad+ \sum_{i\in\mathcal{D}_2} (\widehat{R}_{T,i} - R_{T,i})(R_{Y,i} - D_i\theta_0) + \sum_{i\in\mathcal{D}_2} (\widehat{R}_{T,i} - R_{T,i})(\widehat{R}_{Y,i} - R_{Y,i}).
\end{align*}
Conditional on $\mathcal{D}_1$, the first three terms above are mean zero and hence,
\[
\mathbb{E}[\mathbb{Z}_n(\theta_0)\big|\mathcal{D}_1] = \sum_{i\in\mathcal{D}_2} \mathbb{E}[(\widehat{g}(X_i) - g_0(X_i))(\widehat{m}(X_i) - m_0(X_i))\big|\mathcal{D}_1].
\]
This can be bounded as
\[
|\mathbb{E}[\mathbb{Z}_n(\theta_0)\big|\mathcal{D}_1]| ~\le~ |\mathcal{D}_2|\|\widehat{g} - g_0\|_{2,n}\|\widehat{m} - m_0\|_{2,n},
\]
where for $P_{X_i}(\cdot)$ representing the probability measure of $X_i$,
\begin{align*}
\|\widehat{g} - g_0\|_{2,n}^2 &:= |\mathcal{D}_2|^{-1}\sum_{i\in\mathcal{D}_2} \int (\widehat{g}(x) - g_0(x))^2dP_{X_i}(x),\\
\|\widehat{m} - m_0\|_{2,n}^2 &:= |\mathcal{D}_2|^{-1}\sum_{i\in\mathcal{D}_2} \int (\widehat{m}(x) - m_0(x))^2dP_{X_i}(x).
\end{align*}
Hence if $|\mathcal{D}_2|^{1/2}\|\widehat{g} - g_0\|_{2,n}\|\widehat{m} - m_0\|_{2,n} = o_p(1)$, then inequality~\eqref{eq:median-bias-partial-regression} implies that the median bias of $\widehat{\theta}_n$ converges to zero. Further, if $|\mathcal{D}_2|^{1/2}\|\widehat{g} - g_0\|_{2,n}\|\widehat{m} - m_0\|_{2,n} = O_p(1)$, then the median bias of $\widehat{\theta}_n$ is bounded away from $0$ to $1/2$.
\section{Conclusion}
In this note, we proved that median bias of several M/Z-estimators can be controlled under conditions weaker than those required for convergence in distribution of these estimators. The control of the median bias implies that the recently proposed confidence interval methodology HulC~\citep{Kuchibhotla2020HulC} can be applied to these estimators. Note that without convergence in distribution none of the usual methods of inference, including bootstrap and subsampling, apply. 
\bibliographystyle{apalike}
\bibliography{references}

\begin{thebibliography}{}

\bibitem[Bentkus et~al., 1997]{bentkus1997berry}
Bentkus, V., Bloznelis, M., and G{\"o}tze, F. (1997).
\newblock A {B}erry--{E}ss{\'e}en bound for {M}-estimators.
\newblock {\em Scandinavian journal of statistics}, 24(4):485--502.

\bibitem[Chernozhukov et~al., 2018]{chernozhukov2018double}
Chernozhukov, V., Chetverikov, D., Demirer, M., Duflo, E., Hansen, C., Newey,
  W., and Robins, J. (2018).
\newblock Double/debiased machine learning for treatment and structural
  parameters.

\bibitem[H{\"o}rmann, 2009]{hormann2009berry}
H{\"o}rmann, S. (2009).
\newblock {B}erry-{E}sseen bounds for econometric time series.
\newblock {\em ALEA Lat. Am. J. Probab. Math. Stat}, 6:377--397.

\bibitem[Knight, 1998]{knight1998limiting}
Knight, K. (1998).
\newblock Limiting distributions for ${L}_1$ regression estimators under
  general conditions.
\newblock {\em Annals of statistics}, pages 755--770.

\bibitem[Knight, 1999]{knight1999asymptotics}
Knight, K. (1999).
\newblock Asymptotics for ${L}_1$-estimators of regression parameters under
  heteroscedasticity.
\newblock {\em Canadian Journal of Statistics}, 27(3):497--507.

\bibitem[Kuchibhotla et~al., 2021]{Kuchibhotla2020HulC}
Kuchibhotla, A.~K., Balakrishnan, S., and Wasserman, L. (2021).
\newblock The {HulC}: Confidence regions from convex hulls.
\newblock {\em arXiv preprint arXiv:2105.14577}.

\bibitem[Petrov, 2012]{petrov2012sums}
Petrov, V.~V. (2012).
\newblock {\em Sums of independent random variables}, volume~82.
\newblock Springer Science \& Business Media.

\bibitem[Van~der Vaart, 2000]{van2000asymptotic}
Van~der Vaart, A.~W. (2000).
\newblock {\em Asymptotic statistics}, volume~3.
\newblock Cambridge university press.

\end{thebibliography}
\end{document}